\documentclass[12pt]{amsart}
\usepackage{amssymb,amscd,verbatim}

\swapnumbers
\newtheorem{thm}{Theorem}[subsection]
\newtheorem{propose}[thm]{Proposition}
\newtheorem{lemma}[thm]{Lemma}

\theoremstyle{definition}
\newtheorem{defn}[thm]{Definition}
\newtheorem{remark}[thm]{Remark}
\newtheorem{example}[thm]{Example}

\renewcommand{\d}{\mbox{\LARGE $\cdot $}}
\newcommand{\Xs}{X_{\d}}            
\renewcommand{\hat}{\widehat}
\newcommand{\Lie}{{\rm Lie}\,}       
\newcommand{\Hom}{\operatorname{Hom}}      
\newcommand{\Ext}{\operatorname{Ext}}      
\newcommand{\M}{\mathcal{M}_1}   

\newcommand{\ihom}{{\rm\underline{Hom}}}  

\newcommand{\C}{\mathbb{C}}     
\newcommand{\Q}{\mathbb{Q}}     
\newcommand{\Z}{\mathbb{Z}}     

\newcommand{\G}{\mathbb{G}}     
\newcommand{\HH}{\mathbb{H}}    
\newcommand{\bPic}{{\rm\mathbb{P}ic}}    

\renewcommand{\ker}{\operatorname{Ker}}  
\newcommand{\gr}{\operatorname{gr}}        
\newcommand{\Pic}{{\rm Pic}}     

\newcommand{\by}[1]{\stackrel{#1}{\rightarrow}}
\newcommand{\longby}[1]{\stackrel{#1}{\longrightarrow}}

\newcommand{\df}{\mbox{\,${:=}$}\,}
\newcommand{\ie}{{\it i.e.\/},\ }
\newcommand{\cf}{{\it cf.\/}\ }
\newcommand{\eg}{{\it e.g.\/},\ }
\newcommand{\et} {\mbox{\scriptsize{\rm {\'e}t}}}

\newcommand{\an}{\mbox{\scriptsize{\rm an}}}
\newcommand{\fr}{\mbox{\scriptsize{\rm fr}}}

\renewcommand{\bar}{\overline}
\newcommand{\into}{\hookrightarrow}

\newcommand{\sZ}{\mbox{\scriptsize{$\Z$}}}   
\newcommand{\sC}{\mbox{\scriptsize{$\C$}}}   
\newcommand{\sQ}{\mbox{\scriptsize{$\Q$}}}   
\renewcommand{\c}{\mbox{\scriptsize{$\circ$}}} 

\newcommand{\onto}{\mbox{$\to\!\!\!\!\to$}}

\newcommand{\boxtensor}{\def\boxtimesten{\Box\kern-7.59pt\raise1.2pt
\hbox{$\times$} }}                                  

\newcounter{elno}                   

\newcommand{\cO}{\mathcal{O}}

\begin{document}

\title{Formal Hodge Theory}
\author{ Luca Barbieri-Viale}
\address{Dipartimento di Matematica Pura e Applicata, Universit\`a degli Studi di Padova\\ Via G. Belzoni, 7\\Padova -- I-35131\\ Italy}
\email{barbieri@math.unipd.it}
\date{November 22, 2005} 

\begin{abstract} We introduce {\em formal}\, (mixed) Hodge structures  (of level $\leq 1$) in such a way that the Hodge realization of Deligne's 1-motives extends to a realization from Laumon's 1-motives to formal Hodge structures (of level $\leq 1$) providing an equivalence of categories.
\end{abstract}
 
\maketitle
 
Let ${\rm MHS}_1^{\fr}$ denote the category of torsion free graded polarizable mixed Hodge structures of level $\leq 1$. We have a nice algebraic description of this category {\it via}\, 
$\M^{\fr}$ the category of Deligne's 1-motives  \cite{D} (\cf also \cite{BRS}, including torsion, one obtains 1-motives with torsion describing ${\rm MHS}_1$). Actually, Deligne's Hodge realization provide an equivalence $$T_{Hodge}: \M^{\fr}\longby{\simeq} {\rm MHS}_1^{\fr}$$ such that Cartier duality on $\M^{\fr}$ is transformed in $\ihom (-, \Z (1))$ on ${\rm MHS}_1^{\fr}$.
Moreover, we have a natural generalization of Deligne's 1-motives due to Laumon \cite{LAU}. A Laumon 1-motive $M \df [F \by{u} G]$ is a commutative formal group $F=F^0\times F_{\et}$, with torsion free \'etale part $F_{\et}$, a commutative connected algebraic group $G$ and a map of abelian fppf-sheaves $u: F\to G$. Let $\M^{a,\fr}$ denote the category of Laumon's 1-motives and refer to its objects as 1-motives for short. Note that Cartier duality on $\M^{\fr}$ canonically extends to $\M^{a,\fr}$ (see \cite{LAU}).\\

The purpose of this note is to introduce the abelian category ${\rm FHS}_1$ of {\em formal}\, mixed Hodge structures   (of level $\leq 1$) in order to extend the Hodge realization $T_{Hodge}$ of Deligne's 1-motives $\M^{\fr}$ to a realization $T _{\oint}$ from Laumon's 1-motives $\M^{a,\fr}$ to ${\rm FHS}_1^{\fr}\subset {\rm FHS}_1$. 
We have that ${\rm MHS}_1^{\fr}\subset {\rm FHS}_1^{\fr}$ in a canonical way, \ie there is a fully faithful  embedding such that the natural involution (Cartier duality) on ${\rm MHS}_1^{\fr}$ extends to an involution on ${\rm FHS}_1^{\fr}$.\\

For the sake of exposition we here confine our study to level $\leq 1$ mixed Hodge structures. However, it is conceivable and suitable to consider formal mixed Hodge structures with arbitrary Hodge numbers: generalizing our definition below it's not that difficult (we will treat such a matter nextly, \cf \cite[2.12]{THS} for the general setting). For example, enriched Hodge structures \cite{BSE} (of level $\leq 1$) can easily be recovered as ``special'' formal Hodge structures (see also \cite{BAB} for details). In \cite{BAB} we are also providing a ``sharp'' De Rham realization generalizing Deligne's 
construction of De Rham realization in \cite{D}.
The main result of this paper can be summarized in the following way.

\subsection*{\bf Theorem}{\em There is an equivalence of categories with involution
$$T _{\oint}: \M^{a,\fr}\longby{\simeq} {\rm FHS}_1^{\fr}$$
between Laumon's 1-motives and torsion free formal Hodge structures (of level $\leq 1$)
providing a diagram
$$\begin{array}{ccc}
\M^{\fr}&\by{\simeq}& \mbox{\rm MHS}_1^{\fr}\\
\uparrow\downarrow & &\uparrow\downarrow\\
\M^{a,\fr}&\by{\simeq}& \mbox{\rm FHS}_1^{\fr}
\end{array}$$
where
\begin{itemize}
\item $\M^{\fr}\into \M^{a,\fr}$ and ${\rm MHS}_1^{\fr}\into {\rm FHS}_1^{\fr}$ are canonical inclusions,
\item $\M^{a,\fr}\to \M^{\fr}$ and ${\rm FHS}_1^{\fr}\to {\rm MHS}_1^{\fr}$ are ``forgetful functors'' denoted  $( \ \ )\leadsto (\ \ )_{\et}$, which are left inverses of the inclusions,
\item $T _{\oint}(M) $ coincide with $T_{Hodge}(M)$ if $M = M_{\et}$ and, in general, we have a formula
$$T _{\oint}(M)_{\et} = T_{Hodge}(M_{\et}).$$
\end{itemize}\vspace{0.5cm} }
 
The plan of the paper is the following. In Section 1 we introduce the category ${\rm FHS}_1$.
In Section 2 we construct $T _{\oint}$ proving the theorem.

\section{Formal Hodge Structures}\label{form}  
 
 \subsection{Paradigma} Consider a commutative formal group $H=H^0\times H_{\sZ}$ over $\C$  along with a mixed Hodge structure on the \'etale part $H_{\sZ}$, \ie say $H_{\et}\df (H_{\sZ}, W_*, F^*_{Hodge}) \in {\rm MHS}_1$ for short. For the mixed Hodge structure $H_{\et} \in {\rm MHS}_1$ we here denote $H_{\sZ}$ the finitely generated abelian underlying group, along with the weight filtration $W_{-2}\subseteq W_{-1}$ of $H_{\sQ}\df H_{\sZ}\otimes \Q$ and $F^0_{Hodge}\subseteq H_{\sC}\df H_{\sZ}\otimes \C$ the Hodge filtration.
 We say that $H$ is {\it free}\, if the \'etale part of the formal group  is free, so that: $H_{\sZ}=\Z^r $ and $H^0=\hat{\C}^s$ non-canonically. (Note that here $\hat{\C}$ denotes the connected formal additive group).  For $H$ free we also denote
by  $W_{*}H_{\et}$ and $\gr^W_{*}H_{\et}$  the corresponding objects of $ {\rm MHS}_1$.

 \begin{defn}\label{fhsdef}
Define a {\it formal Hodge structure\/} (of level $\leq 1$) as follows: {\it (i)}\, a formal group $H$  such that $H_{\et}\in {\rm MHS}_1$,  {\it (ii)}\, a finite dimensional $\C$-vector space $V$ with a two steps filtration $V^0\subseteq V^1\subseteq  V $ by sub-spaces,  {\it (iii)}\,  a group homomorphism $v: H\to V$ and {\it (iv)}\, a $\C$-isomorphism $\sigma : H_{\sC}/F^0_{Hodge} \by{\simeq} V/V^0$ restricting to an isomorphism
$ W_{-2}H_{\sC}\cong V^1/V^0$. We further assume that the following condition holds: if $v_{\sZ}: H_{\sZ} \to V$ is the induced map, $c : H_{\sZ} \to H_{\sC}/F^0_{Hodge}$ is the canonical map and ${\rm pr} : V\onto V/V^0$ is the projection then the following
\begin{equation}\label{cond}
\begin{array}{ccc}
H_{\sZ}&\longby{v _{\mbox{\tiny $\Z$}}}&V\\
\mbox{\tiny c}\downarrow & &\downarrow \mbox{\it\tiny pr}\\
H_{\sC}/F^0_{Hodge}&\longby{\sigma}& V/V^0
\end{array}
\end{equation}
commutes. Denote $(H, V)$ for short such a structure. 
\end{defn}

Define a morphism $\phi$ between $(H, V)$ and $(H', V')$ as follows. We let $\phi \df (f, g)$ be a pair of maps in the following commutative square
\begin{equation}\label{map}\begin{array}{ccc}
H&\longby{v }&V\\
\mbox{\it\tiny f}\downarrow & &\downarrow\mbox{\it\tiny g}\\
H'&\longby{v'}& V'
\end{array}\end{equation}
where $f : H\to H'$ is a homomorphism of formal groups such that $f_{\et} : H_{\et}\to H'_{\et}$ is a map in ${\rm MHS}_1$
and $g :V\to V'$ is a $\C$-homomorphism compatible with the filtrations, \ie $g (V^i)\subseteq V'^{i}$ for $i= 0, 1$.  We further assume that the following diagram commutes 
\begin{equation}\label{condmap}
\begin{array}{ccc}
H_{\sC}/F^0_{Hodge}&\longby{\sigma}& V/V^0\\
\mbox{\tiny $\bar f $}\downarrow & &\downarrow \mbox{\tiny $\bar g$}\\
H_{\sC}'/F^0_{Hodge}&\longby{\sigma'}& V'/V'^0
\end{array}
\end{equation}
where  $\bar f $ and  $\bar g$ are the canonically induced maps.

\begin{defn} Let ${\rm FHS}_1$ denote the {\it category}\, whose objects are $(H, V)$, the morphisms are $\phi = (f, g)$ as above and the composition is given by gluing the squares \eqref{map} (the condition \eqref{condmap} is preserved by gluing). Let ${\rm FHS}_1^{\fr}\subset {\rm FHS}_1$ denote the full subcategory given by $(H, V)$ such that $H$ is free.
\end{defn}
\begin{propose} The category ${\rm FHS}_1$ is abelian. 
A  short  exact   sequence
$$0\to (H, V) \to (H', V') \to (H'',V'')\to 0$$
is given by an exact sequence on each component (formal groups and filtered vector spaces) so that
$$0\to H_{\et} \to H'_{\et}\to H''_{\et}\to 0$$
is exact in ${\rm MHS}_1$.
\end{propose}
\begin{proof} Straightforward. \end{proof} 

\subsection{\'Etale structures} We can recover mixed Hodge structures as follows.
\begin{defn} Define $(H, V)_{\et}\df (H_{\sZ}, V/V^0)$ where $(H_{\sZ})_{\et}\df H_{\et}$, $v_{\et} : H_{\sZ}\to V/V^0$ is the composition of pr and $v_{\sZ}$ (\cf \ref{fhsdef}) and $(V/V^0)^0\df 0\subseteq (V/V^0)^1\df V^1/V^0\subseteq V/V^0$.
Say that a formal Hodge structure is {\it \'etale}\, if $(H, V) = (H, V)_{\et}$, \ie if $H^0=V^0=0$.
\end{defn}
Given $ (H_{\sZ}, W_*, F^*_{Hodge}) \in  {\rm MHS}_1$ there is a natural way to provide an \'etale one as follows. Set $H \df H_{\sZ}$, $H_{\et} \df (H_{\sZ}, W_*, F^*_{Hodge}) $,  $H^0=0$, $V\df H_{\sC}/F^0_{Hodge}$, $V^1 \df W_{-2}H_{\sC}$, $V^0\df 0$, $\sigma$ is the identity and  the map $v\df c$ is induced by the canonical map $t: H_{\sZ} \to H_{\sC}$. Denote 
$$ c(H_{\sZ}, W_*, F^*_{Hodge}) \df (H_{\sZ}, H_{\sC}/F^0_{Hodge})$$
the {\it canonical}\, \'etale formal Hodge structure associated to a mixed Hodge structure, providing a functor $c :   {\rm MHS}_1\to   {\rm FHS}_1$.
\begin{lemma}  \label{fet} The full subcategory ${\rm FHS}_1^{\et}$ of \'etale structures is equivalent to ${\rm MHS}_1$ \emph{via} $c$ and the forgetful functor $(H, V)\mapsto H_{\et}$. The functor $e: (H, V)\mapsto (H, V)_{\et}$ is a left inverse of the  inclusion ${\rm FHS}_1^{\et}\subset {\rm FHS}_1$ and, for $ (H',V') \in {\rm FHS}_1^{\et}$, we have
$$ \Hom ( (H,V), (H',V')) \subseteq \Hom ( (H,V)_{\et}, (H',V') ) $$
where the equality holds if $v (H^0)\subseteq V^0$ (\cf 1.3.1 below).
\end{lemma}
\begin{proof} Actually, for the equivalence, we are easily left to show that if $(H, V)$ is \'etale then $c(H_{\et})\df (H_{\sZ}, H_{\sC}/F^0_{Hodge})\cong (H, V)$. The claimed isomorphism is $(1, \sigma)$ granted by \eqref{cond} since $V^0=H^0=0$.

For the other claims, let $(H, V)\in {\rm FHS}_1$ and  $ (H',V') \in {\rm FHS}_1^{\et}$ and  consider a map $\phi =(f,g): (H,V)\to (H',V') $ whence induced maps $\bar f$ and $\bar g$ and a diagram
$$\begin{array}{ccc}
H &\longby{v }&V\\
\uparrow\downarrow & &\downarrow\\
H_{\sZ}&\longby{v _{\et}}&V/V^0\\
\mbox{\tiny $\bar f$}\downarrow & &\downarrow\mbox{\tiny $\bar g$}\\
H'&\longby{v'}& V'
\end{array}$$
In fact $H'=H'_{\et}$ is \'etale thus $f (H^0) =0$ and $f$ factors through $H_{\sZ}$ yielding $\bar f$ and, similarly, we get a filtered map $\bar g: V/V^0\to V'$ since $V'^0=0$ and $g (V^0) =0$. Now $\bar \phi \df (\bar f, \bar g)$ yields a map by diagram chase. Note that 
if $v (H^0)\subseteq V^0$ then $ (H,V)\to (H,V)_{\et}$ (\cf \eqref{fmncext}) and we can lift back, by composition, any morphism $\phi': (H,V)_{\et}\to (H',V') $ as the condition \eqref{condmap} is tautological.
\end{proof}

\begin{remark} Note that under the equivalence we then get a canonical inclusion $c : {\rm MHS}_1^{\fr}\into {\rm FHS}_1^{\fr}$ such that 
$e : {\rm FHS}_1^{\fr}\to {\rm MHS}_1^{\fr}$ is a left inverse and for $H'\in {\rm MHS}_1^{\fr}$
$$\Hom_{{\rm FHS}_1^{\fr}}( (H,V), (H'_{\sZ},H_{\sC}'/F^0_{Hodge}) ) \subseteq \Hom_{{\rm MHS}_1^{\fr}} ( H_{\et}, H')$$
\end{remark}
 
 \subsection{Connnected structures}  A $\C$-vector space $V$ will be regarded as an object $(0,V)$ of ${\rm FHS}_1$ filtered as $V=V^1=V^0$. Similarly, a formal group $H$ is regarded as an object $(H, 0)$ of ${\rm FHS}_1$ so that $H =H^0\times H_{\sZ}$ and  $H_{\sZ}$ is pure of weight zero. 

For $(H, V)\in {\rm FHS}_1$ we have that $V^0$ is a substructure of $(H, V)$ and we can consider the quotient $(H, V)/V^0= (H, V/V^0)$ in ${\rm FHS}_1$. 
We can also regard $(H, V)_{\et}$ as a substructure of $(H, V)/V^0$ and we obtain a canonical exact sequence
\begin{equation}\label{fmext} 0\to (H, V)_{\et}\to (H, V)/V^0\to H^0\to 0
\end{equation} 

 \begin{defn} Say that $(H,V)\in {\rm FHS}_1$ is {\it connected}\, if $H =H^0$ is connected, \ie if $(H,V)_{\et}=0$. 
 Denote $\pi (H,V) \df (H^0,V)$ the connected structure given by $V=V^1=V^0$ and the restriction of $v$ to $H^0\subseteq H$. Let ${\rm FHS}_1^0$ denote the full subcategory of $ {\rm FHS}_1$ determined by connected structures.
 
 Say that $(H,V)\in {\rm FHS}_1$ is {\it special}\, if $v (H^0)\subseteq V^0$, \ie if $v: H\to V$  restricts to $v^0:H^0\to V^0$.  Denote  ${\rm FHS}_1^s$ the full subcategory of special structures and $ (H,V) ^0 \df (H^0,V^0)\in {\rm FHS}_1^0$ the connected structure determined by $(H,V)\in {\rm FHS}_1^s$. \end{defn} 
 \begin{lemma}  \label{cet}  The functor $ (H, V)\mapsto \pi (H, V)$ is a left inverse of the  inclusion $\iota : {\rm FHS}_1^{0}\subset {\rm FHS}_1$. The category ${\rm FHS}_1^{0}$ is equivalent to the category of linear mappings between finite dimensional $\C$-vector spaces. For $ (H',V') \in {\rm FHS}_1^{0}$ and $(H,V)\in {\rm FHS}_1^s$
$$\Hom ( (H',V'), (H,V)) \cong \Hom ( (H',V'), (H,V)^0)$$
\end{lemma}

\begin{proof} The first claim is clear.  Moreover, the equivalence is provided by $(H, V)\mapsto \Lie (H) \to V$. 
Finally, a map from $(H', V')$ connected to $ (H,V) $ special is given by a commutative square
$$\begin{array}{ccc}
H'&\longby{}&V'\\
\mbox{\tiny $f$}\downarrow\ \ & &\ \ \downarrow\mbox{\tiny $g$}\\
H&\longby{}& V
\end{array}$$
such that $f (H') \subseteq H^0$ and $g (V')\subseteq V^0$.
\end{proof}
\begin{remark} Note that $(H,V)$ with $H_{\et}$ pure of weight zero exists if and only if $V=V^1=V^0$. Thus if $(H, V)$ is special then $(H,V)^0$ is the largest connected formal substructure of $(H, V)$ and we have a {\it non canonical}\, extension
\begin{equation}\label{fmncext} 0\to (H^0, V^0)\to (H,V) \to (H,V)_{\et}\to 0
\end{equation} 
From lemmas~\ref{fet} and \ref{cet} it follows that the functors $ (H, V)\mapsto (H, V)^0$ and $(H,V)\mapsto  (H,V)_{\et}$ are, respectively, a right adjoint 
of ${\rm FHS}_1^{0}\subset {\rm FHS}_1^s$ and a left adjoint
of  ${\rm FHS}_1^{\et}\subset {\rm FHS}_1^s$. However, special structures do have   disadvantages, see 2.2.5 and 2.3.2.
\end{remark}
 \begin{propose} The category ${\rm FHS}_1^0$ forms a Serre abelian subcategory of ${\rm FHS}_1$ yielding the extension
 $$0\to {\rm FHS}_1^0\by{\iota}  {\rm FHS}_1 \by{e}{\rm MHS}_1\to 0$$
where $\pi\iota =1$ and $e c =1$.
 \end{propose}
 \begin{proof} It follows from the lemmas~\ref{fet},  \ref{cet} and \eqref{fmext}. In fact, it is clear (\cf 1.1.3) that ${\rm FHS}_1^0$ forms a Serre subcategory. Since $e ({\rm FHS}_1^0)=0$ we have a factorisation 
 $\bar e : {\rm FHS}_1/{\rm FHS}_1^0\to {\rm MHS}_1$ {\it via}\, the canonical projection $t : {\rm FHS}_1 \to {\rm FHS}_1/{\rm FHS}_1^0$ and the equivalence 
 ${\rm FHS}_1^{\et}\cong {\rm MHS}_1 $. Since $e = \bar e t$ and $ec =1$ then $\bar e tc = 1$. We also have $t c \bar e \cong 1$ since applying $t$ to 
 \eqref{fmext} for $(H, V)\in  {\rm FHS}_1$ we get a natural isomorphism
 $$t c(H_{\et})\cong t (H,V)_{\et}\cong t(H, V)$$ \end{proof}
 
 \subsection{Construction} We provide a Laumon 1-motive out of a {\it free}\, formal mixed Hodge structure (of level $\leq 1$). The construction is similar to \cite[p. 55-56]{D}. 
 
For $(H, V)\in {\rm FHS}_1^{\fr}$ the Laumon 1-motive $\overrightarrow{(H, V)} \df [F\by{u}G]$ functorially associated to  $(H, V)$ is given as follows. Set $F\df H^0\times \gr^W_0(H_{\sZ})$. Since \eqref{cond} holds true $W_{-1}(H_{\sZ})$ injects in $V$ {\it via}\, $v_{\sZ}: H_{\sZ}\to V$ in such a way  that $W_{-1}(H_{\sZ})\cap V^0 =0$. Set $G (\C)\df V/W_{-1}(H_{\sZ})$ obtaining a diagram
\begin{equation}\label{deflaum}
\begin{CD}
0@>>> W_{-1}(H_{\sZ}) @>{}>> H @>{}>>F@>>> 0\\ 
&&@V{||}VV @V{v}VV@V{u}VV\\
0@>>>W_{-1}(H_{\sZ})@>>> V @>>> G (\C)@>>> 0
\end{CD}
\end{equation} 
where $u$ is just induced by $v$. Regarding the complex group $G (\C)$ we then have it in
a diagram 
\[\begin{CD}
&& && 0 &&
0&&\\ 
&& && @V{}VV@V{}VV\\
&& && V^0 @>{=}>>
V^0 &&\\ 
&& && @V{}VV@V{}VV\\
0@>>> W_{-1}(H_{\sZ}) @>{\tiny v_{\Z}}>> V @>{}>>G (\C)@>>> 0\\ 
&&@V{||}VV @V{}VV@V{}VV\\
0@>>>W_{-1}(H_{\sZ})@>{c}>> H_{\C}/F^0_{Hodge}@>>> J (W_{-1}(H_{\et}))@>>> 0\\
&& && @V{}VV@V{}VV\\
&& && 0 &&
0&&\\ 
\end{CD}\]
obtained {\it via}\, $\sigma$ and \eqref{cond}. This is showing that $G (\C)$ is an extension of the complex torus $J (W_{-1}(H_{\et}))$ by a $\C$-vector group. Thus, by G.A.G.A., we get the algebraic group $G$. 
 
\section{Formal Hodge realization}
\subsection{Paradigma} For a Laumon 1-motive $M = [F\by{u}G]\in \M^{a,\fr}$ over a field $k$ (algebraically closed of characteristic zero) we here denote $F = F^0\times F_{\et}$ the formal group where $F_{\et}$ is further assumed torsion free. Denote $V (G)\df \G_a^n\subseteq G$ the additive factor and display the connected algebraic group $G$ as an extension
\begin{equation}\label{gmext} 0\to V (G)\to G \to G_{\times}\to 0
\end{equation} where $G_{\times}$ is the semi-abelian quotient.  The algebraic group $G_{\times}$ is an extension of an abelian  variety $A$ by a torus $T$.
 \begin{defn} For $M= [F\by{u}G]\in \M^{a,\fr}$ set $M_{\et}\df [F_{\et}\by{u_{\et}} G_{\times}]\in \M^{\fr}$. Say that  $M$ is {\it \'etale}\, if $M=M_{\et}$, \ie it is a Deligne 1-motive. Say that $M$ is {\it connected}\, if $M_{\et}=0$, \ie $F = F^0$ is connected and $G = V(G)$ is a vector group. Say that $M$ is {\it special}\, if 
 $u (F^0)\subseteq V(G)$ and set $M^0\df [F^0\to V(G)]$. \end{defn} 
\begin{lemma}\label{etform}
The functor $M\mapsto M_{\et}$ is a left inverse of the inclusion $\M^{\fr} \subset \M^{a,\fr}$ of Deligne's 1-motives and for $M'\in \M^{\fr} $ we have
$$\Hom  (M, M') \subseteq \Hom (M_{\et}, M') $$
  If $M$ is special we then get an extension
 \begin{equation}\label{ncformext} 0\to M^0\to M \to M_{\et}\to 0
\end{equation} 
such that if $M'$ is \'etale then 
$\Hom (M_{\et}, M') \cong \Hom  (M, M') $
and if $M'$ is connected then
$\Hom (M',M^0) \cong \Hom  (M', M) $.
\end{lemma}
\begin{proof} Let $M = [F\by{u} G]\in \M^{a,\fr}$, $M' = [F'\by{u'} G']\in \M^{\fr}$. Let $(f, g):M\to M'$ be a map. Then get a diagram (\cf the proof of \ref{fet})
$$\begin{array}{ccc}
F&\longby{u }&G\\
\uparrow\downarrow & &\downarrow\\
F_{\et}&\longby{u_{\et}}&G_{\times}\\
\mbox{\tiny $\bar f$}\downarrow & &\downarrow\mbox{\tiny $\bar g$}\\
F'&\longby{u'}& G'
\end{array}$$
where $\bar f$ and $\bar g$ are the induced maps since $M'$ is \'etale, yielding a map
$(\bar f, \bar g):M_{\et}\to M'$. In fact, $\Hom (F, F') = \Hom (F_{\et}, F')$ because $F'$ is \'etale and $F^0$ is mapped to zero and $\Hom (G, G') = \Hom (G_{\times}, G')$ because $\Hom (\G_a,\G_m)=\Hom (\G_a, A)=0$ and $G'$ is semi-abelian. Moreover, $M\to M_{\et}$ if $M$ is special, yielding  \eqref{ncformext}. For the isomorphisms then note that $\Hom (M^0, M')=0$ if $M'$ is \'etale and, equivalently, $\Hom (M', M_{\et})=0$ if $M'$ is connected.
\end{proof}

In general,  we can regard $M_{\et}$ as a sub-1-motive of $M/V(G)$ and we obtain (\cf \eqref{fmext}) a canonical exact sequence
\begin{equation}\label{lmext} 0\to M_{\et}\to M/V(G) \to F^0[1]\to 0
\end{equation} 
Denote $M_{\et}^{\natural}=[F_{\et}\by{u^{\natural}}G^{\natural}]\in \M^{a,\fr}$ (\cf \cite{D})  the universal $\G_a$-extension of $M_{\et}$. 
The algebraic group $ G^{\natural}$ can be represented by an extension 
\begin{equation}\label{unext}
0\to \Ext (M_{\et},\G_a)^{\vee}\to G^{\natural}\to G_{\times}\to 0
\end{equation}
where $\Ext (M_{\et},\G_a)^{\vee}$ is given by the dual vector space of $\G_a$-extensions of $M_{\et}$. The map $u^{\natural}: F_{\et}\to G^{\natural}$ is a canonical lifting of $u_{\et} : F_{\et}\to G_{\times}$.\\

Set $k =\C$. Recall that Deligne's Hodge realization  (see \cite{D})  $$T_{Hodge }(M_{\et})\df (H_{\sZ },W_*, F^0_{Hodge})$$  of $M_{\et}$ is given by the pull-back
$$\begin{array}{ccc}
H_{\sZ }&\longby{\bar v _{\tiny \Z}}&\Lie (G_{\times})\\
\downarrow & &\ \ \downarrow\mbox{\tiny $\exp$}\\
F_{\et}&\longby{u_{\et}}& G_{\times}
\end{array}$$
Here $W_{-1} \df H_1(G_{\times})$, $W_{-2} \df H_1(T)$ and 
$$F^0_{Hodge}\df \ker (H_{\sC} \to \Lie (G_{\times}))$$
\begin{lemma}\protect{\rm (\cite[10.1]{D})}\label{derham}
For $k= \C$ we have an isomorphism 
$$M_{\et}^{\natural}\cong [H_{\sZ}/W_{-1}\by{\bar t} H_{\sC}/W_{-1}]$$ here $\bar t$ is the induced  map $t : H_{\sZ}\to H_{\sC}\mod W_{-1}(H_{\sZ})$. \end{lemma}

Actually (see \cite[10.1.8]{D}) we have a bifiltered isomorphism (\ie ``periods'')
$$\tau :  \Lie (G^{\natural})\by{\simeq}H_{\sC }$$ such that 
\begin{equation}\label{periods}
\begin{array}{ccccc}
H_{\sZ }&\longby{v ^{\natural}}&\Lie (G^{\natural})& \by{\tau}&H_{\sC}\\
\scriptstyle{||}& &\ \ \downarrow\scriptstyle{}& &\downarrow\\
H_{\sZ }&\longby{\bar v _{\tiny \Z}}&  \Lie (G_{\times})& \by{\bar \tau}&H_{\sC}/F^0_{Hodge}
\end{array}
\end{equation}
commutes. Here $t = \tau v^{\natural}$ where $v^{\natural}$ is the canonical  map induced by $u^{\natural}$, $\Lie (G^{\natural})$ is the pullback of \eqref{unext} along exp, $H_1(G^{\natural})\cong H_1(G_{\times})= W_{-1}(H_{\sZ})$ and $\Ext (M_{\et},\G_a)^{\vee}\cong F^0_{Hodge}$. 

\begin{example} \protect{\rm (\cf \cite[1.1 \& 3.3]{THS})} For $X$ proper over a field $k$, char $(k)=0$, set $G\df \Pic_{X/k}^0$ and let $M = [0\to G]$ be the  corresponding 1-motive. Here $G_{\times}\cong \bPic_{\Xs /k}^0$ and $G^{\natural}\cong \bPic_{\Xs /k}^{\natural, 0}$ are given by simplicial Pic and $\natural-\Pic$ functors of a smooth proper hypercovering $\Xs$ of $X$. Thus $H_{\sZ }= H^1(X_{\rm an},\Z )$, $\Lie (G^{\natural})= H^1_{DR}(X)$ and $\tau : H^1_{DR}(X)\cong H^1(X_{\rm an},\C )$ by cohomological descent over $k =\C $.
\end{example}

\subsection{Construction} Extending Deligne's Hodge realization for a given Laumon 1-motive $M = [F \by{u} G]$ over $\C$ consider the pull-back $T_{\oint}(F)$ of $u: F \to G$ along $\exp :\Lie (G)\to G$, \ie 
$$\begin{array}{ccc}
T_{\oint}(F) &\longby{v}&\Lie (G)\\
\downarrow & &\ \ \ \downarrow\mbox{\tiny $\exp$}\\
F&\longby{u}& G
\end{array}$$
Here $T_{\oint}(F)$  is a formal group and we get  a natural group homomorphism $v: T _{\oint}(F)\to \Lie (G)$. We are going to show that $$T _{\oint}(M)\df (T _{\oint}(F),\Lie (G))\in {\rm FHS}_1^{\fr}$$ is a formal Hodge structure. 
Note that if $M$ is connected then $T _{\oint}(M) = M$.
\begin{remark} The additional data coming from $\Lie$ is really needed if  we allow additive factors! For example, let $W\to V$ be a linear map between $\C$-vector spaces, and let $M = [\hat{W}\by{u} V]$ be the  induced 1-motive where $\hat{W}$ is the
formal completion at the origin  (\cf \cite[5.2.5]{LAU}).  Note that all connected 1-motives are obtained in this way (see \ref{cet}).
For any embedding $V\subsetneq V^{\prime}$ of vector spaces, we
obtain another 1-motive $M^{\prime} = [\hat{W}\by{u'}  V^{\prime}]$ such that $M\subsetneq M^{\prime}$. For both $M$ and $M'$
then $T_{\oint}(\hat{W})$ is the infinitesimal group $\hat{W}$, $\ker (u)=\ker (u')$ and we cannot distinguish $M$ by $M'$ out of  the formal group only. 
\end{remark}
\begin{lemma} \label{fr1}
We have that $T_{\oint}(F)$ is the formal group $F^0\times H_{\sZ}$ such that $H_{\sZ}$ is the above extension of $F_{\et}$ by $H_1(G_{\times})$.\end{lemma} 
\begin{proof} Since formal groups are closed under extensions (\cf  \cite[4.3.1]{LAU}) $T_{\oint}(F)$ is a formal group, \ie it is, by construction, an extension of $F$ by $H_1(G)$. Observe that \eqref{gmext} yields $\Lie (G)$ as the pullback of $\Lie (G_{\times})$ along $\exp$ and 
$H_1(G)\cong H_1(G_{\times})$. We then get a natural identification of $H_{\sZ}$ with the \'etale part of $T _{\oint}(F)$, \ie with the pullback of $F_{\et}\into F$ along $T _{\oint}(F)\to F$.
\end{proof}
\begin{lemma} \label{fr2}
 If  $\sigma \df \bar{\tau}^{-1} :  H_{\sC}/F^0_{Hodge}\by{\simeq} \Lie(G_{\times})$ is the isomorphism induced by \eqref{periods} then $\sigma$ restricts to $W_{-2}(H_{\sC})\cong \Lie (T)$ and  the following
 $$
\begin{array}{ccc}
H_{\sZ}&\longby{v _{\mbox{\tiny $\Z$}}}&\Lie (G)\\
\mbox{\tiny c}\downarrow & &\downarrow \mbox{\it\tiny pr}\\
H_{\sC}/F^0_{Hodge}&\longby{\sigma}& \Lie (G_{\times})
\end{array}
$$
commutes (here $v_{\sZ}$ is the restriction of $v$ and $c$ is the canonical map \cf \eqref{cond}).
\end{lemma} 
\begin{proof} Note that $c = \bar \tau \c \bar v_{\sZ}$ in \eqref{periods} and  $\bar v_{\sZ} = {\rm pr} \c v _{\sZ}$ by Lemma~\ref{fr1}.
\end{proof}

 \begin{defn} Denote $T _{\oint}(M)$ the formal Hodge structure $ (H,V)\in {\rm FHS}_1^{\fr}$ where {\it (i)}\,  $H \df T_{\oint}(F) = F^0\times H_{\sZ}$, $H_{\et}\df T_{Hodge}(M_{\et})$, granted by Lemma~\ref{fr1}, {\it (ii)}\,  $V\df \Lie (G)$, $V^1\df \Lie (T) + V(G)$ and $V^0\df V (G)$, {\it (iii)}\, the map $v: T_{\oint}(F) \to \Lie (G)$ defined above, and {\it (iv)}\, the isomorphism $\sigma \df \bar{\tau}^{-1}  : H_{\sC}/F^0_{Hodge}\by{\simeq} \Lie(G_{\times})$ providing \eqref{cond} by Lemma~\ref{fr2}.
 \end{defn} 
We then have 
$T_{\oint}(M)_{\et} = T_{Hodge}(M_{\et})\in {\rm MHS}_1^{\fr}$
and the construction is clearly functorial (since the diagram \eqref{periods} is natural) providing a functor 
$$T_{\oint}:  \M^{a,\fr}\longby{ }\mbox{\rm FHS}_1^{\fr}$$
such that $T_{\oint}(M) = T_{Hodge}(M)$ if $M$ is \'etale ({\it via}\, \ref{fet}) and $T_{\oint}(M) = M$ if $M$ is connected.
\begin{remark} Note that by applying $T_{\oint}$ to \eqref{lmext} we get \eqref{fmext}, the extension \eqref{ncformext} yields \eqref{fmncext} and $M$ is special $\iff $  $T_{\oint}(M)$ is special.
\end{remark}

\subsection{Conclusion} Summarizing up, see also 1.2 and 1.4, the theorem is proven, \eg in order to show that $T_{\oint}$ yields an equivalence of categories we can argue as in \cite[10.1.3]{D}. For $(H,V)\in {\rm FHS}_1^{\fr}$ we have constructed, in 1.4, a 1-motive $$\overrightarrow{(H,V)}\df [H^0\times \gr_0^W (H_{\sZ})\to V/W_{-1}(H_{\sZ})]$$
It is clear that $T_{\oint} (\overrightarrow{(H,V)}) \cong (H, V) $, see \eqref{deflaum}, which is natural in $(H, V)$. Conversely, for $M = [F\to G]$ we have 
$T_{\oint} (M) \df (T_{\oint}(F), \Lie (G) )$ such that $\overrightarrow{T_{\oint} (M)} \cong M$ functorially in $M$ by construction. One obtains a duality on ${\rm FHS}_1^{\fr}$ after Cartier duality on $\M^{a,\fr}$  by defining
$$T_{\oint} (M)^{\vee}\df T_{\oint} (M^{\vee})$$
The lemmas \ref{fet} and \ref{etform} further explain the diagram of the main theorem and the remaining claims.

\begin{example} \protect{\rm (\cf 2.1.4)} For $X$ proper over $\C$ and 
$M = [0\to \Pic_{X/\sC}^0]$ we have $T_{\oint} (M) = (H^1(X_{\an },\Z (1)),H^1(X,\cO_X))$. Here we have  $M_{\et} = [0\to \bPic_{\Xs/\sC}^{0}]$ and a projection
$$\begin{array}{ccc}
\Lie \Pic_{X/\sC}^0&\longby{\simeq}& H^1(X,\cO_X)\\
\downarrow & &\downarrow\\
\Lie \bPic_{\Xs/\sC}^{0}&\longby{\simeq}& \HH^1(\Xs,\cO_{\Xs})
\end{array}
$$
with kernel the additive factor of $\Pic_{X/\sC}^0$.
Further considering $M^{\natural}_{\et} = [0\to \bPic_{\Xs/\sC}^{\natural, 0}]$ and $T_{\oint} (M^{\natural}_{\et}) = (H^1(X_{\an },\Z (1)),H^1_{DR}(X))$ we get  the extension
$$0\to F^0_{Hodge}\to T_{\oint} (M^{\natural}_{\et}) \to  T_{Hodge} (M_{\et}) \to 0$$
 \end{example}
\begin{remark}
Note that in 2.3.1 $M$ is special but the dual $M^{\vee}$ is not special! Another more striking example is given by taking an abelian variety $X$ and looking at the special 1-motive $[0\to \Pic_{X/\sC}^{\natural, 0}]$ which is the universal extension of the dual of $X$. The Cartier dual 
$$[0\to \Pic_{X/\sC}^{\natural, 0}]^{\vee} = [\hat{X}\to X]$$
is not special. Actually, in general, the Cartier dual of a connected 1-motive is connected  and the dual of \'etale is \'etale but the Cartier dual of $M$ special just fits in an extension 
$$0\to M^{\vee}_{\et}\to M^{\vee} \to (M^0)^{\vee}\to 0$$ 
dual to \eqref{ncformext}.\end{remark}

 \end{document}